\documentclass[12pt]{article}
\usepackage{amssymb,amsfonts,amsmath, psfrag,eepic,colordvi,epsfig}
by \topskip 
\numberwithin{equation}{section}
\newtheorem{theorem}{Theorem}[section]

\newtheorem{lemma}[theorem]{Lemma}

\begin{document}
\parskip 6pt

\pagenumbering{arabic}
\def\sof{\hfill\rule{2mm}{2mm}}
\def\ls{\leq}
\def\gs{\geq}
\def\SS{\mathcal S}
\def\qq{{\bold q}}
\def\MM{\mathcal M}
\def\TT{\mathcal T}
\def\EE{\mathcal E}
\def\lsp{\mbox{lsp}}
\def\rsp{\mbox{rsp}}
\def\pf{\noindent {\it Proof.} }
\def\mp{\mbox{pyramid}}
\def\mb{\mbox{block}}
\def\mc{\mbox{cross}}
\def\qed{\hfill \rule{4pt}{7pt}}
\def\block{\hfill \rule{5pt}{5pt}}

\begin{center}
{\Large\bf Matrix Identities on

 Weighted Partial Motzkin Paths} \vskip
6mm
\end{center}

\begin{center}
\small William Y.C. Chen$^1$,   Nelson Y. Li$^2$, Louis W. Shapiro$^3$ and Sherry H. F. Yan$^4$\\
[2mm] \small $^{1,2,4}$Center for Combinatorics, LPMC, Nankai University, Tianjin 300071, P.R. China \\
[2mm] \small $^3$Department of Mathematics,
                 Howard University, Washington, DC 20059, USA \\
[2mm] \small  $^1$chen@nankai.edu.cn,  $^2$nelsonli@eyou.com, $^3$lshapiro@howard.edu, $^4$huifangyan@eyou.com \\
\end{center}

\vskip10mm

\noindent {\bf Abstract.}  We give a  combinatorial interpretation
of a matrix identity on Catalan numbers and the sequence $(1, 4,
4^2, 4^3, \ldots)$ which has been derived by  Shapiro, Woan and
Getu by using Riordan arrays.  By giving a bijection between
weighted partial Motzkin paths with an elevation line and weighted
free Motzkin paths, we find a matrix identity on the number of
weighted Motzkin paths and the sequence $(1, k, k^2, k^3, \ldots)$
for any $k \geq 2$. By extending this argument to partial Motzkin
paths with multiple elevation lines, we give a combinatorial proof
of an identity recently obtained by Cameron and Nkwanta. A matrix
identity on colored Dyck paths is also given, leading to a matrix
identity for the sequence $(1, t^2+t, (t^2+t)^2, \ldots)$.



\noindent {\sc Corresponding Author:} William Y. C. Chen,
chen@nankai.edu.cn


\vskip10mm
\section{Introduction}

This paper is motivated the following matrix identity obtained by
 Shapiro, Woan
and Getu \cite{shapirotc} in their study of the moments of a
Catalan triangle \cite{chapman, shapiro,sulanke}:
\begin{equation}\label{eq2.2}
\begin{bmatrix} 1 \\2 & 1\\5 & 4 & 1\\14 & 14 & 6 & 1\\42 & 48 & 27 & 8 & 1\\& & \cdots &&&\ddots  \end{bmatrix}
 \begin{bmatrix} 1 \\ 2 \\ 3 \\ 4 \\ 5 \\ \vdots
\end{bmatrix} = \begin{bmatrix} 1 \\ 4 \\ 4^2 \\ 4^3 \\ 4^4 \\ \vdots
\end{bmatrix},
\end{equation}
where  the first column of the first matrix is the Catalan number
$C_n={1\over n+1} {2n \choose n}$ and $a_{i,j}$ (the entry in the
$i$-th row and $j$-th column) is determined by the following
recurrence relation for $j\geq 2 $:
\begin{equation} \label{r2.2}
a_{i,j}=a_{i-1, j-1}+2a_{i-1,j}+a_{i-1,j+1}   .
\end{equation}
Another proof of the above identity is given by Woan, Shapiro and
Rogers \cite{woan} while computing the areas of
parallelo-polyominos via generating functions.

The first result of this paper is a combinatorial interpretation
of the identity (\ref{eq2.2}) in terms of Dyck paths.

One main objective of this paper is to give a matrix identity that
extends the sequence $(1, 4, 4^2, 4^3, \ldots)$ to $(1, k, k^2,
k^3, \ldots)$ in (\ref{eq2.2}). The following matrix identity was
proved by
 Cameron and Nkwanta \cite{cn} that arose in a study of elements of order $2$ in
  {\em Riordan groups} \cite{aigner, shapiroba, shapirotr,
sprugnoli}:
\begin{equation}\label{eq1.1}
\begin{bmatrix} 1 \\3 & 1\\11 & 6 & 1\\45 & 31 & 9 & 1\\197 & 156 & 60 & 12 & 1 & &\\& & \cdots &&&\ddots  \end{bmatrix}
\begin{bmatrix} 1 \\ 3 \\ 7 \\ 15 \\ 31 \\ \vdots
\end{bmatrix} = \begin{bmatrix} 1 \\ 6 \\ 6^2 \\ 6^3 \\ 6^4 \\
\vdots
\end{bmatrix},
\end{equation}
where the entry $a_{i,j}$ ($i$th row and $j$th column) in the
above matrix satisfies the recurrence relation
\begin{equation} \label{r1.1}
a_{i,j}=a_{i-1,j-1}+3a_{i-1,j}+2a_{i-1,j+1}
\end{equation}
for $j\geq 2$ and the $a_{i,1}$ is the $i$-th {\em little
Schr\"{o}der number} $s_i$ (sequence A001003 in \cite{sloane}),
which counts Schr\"{o}der paths of length $2(i+1)$. A {\em
Schr\"oder path} is a lattice path starting at (0, 0) and ending
at $(2n, 0)$ and using steps $H=(2, 0)$, $U=(1, 1)$ and $D=(1,
-1)$
 such that no steps are below the $x$-axis and there
 are no peaks at level one. Imposing this last peak condition
 gives us little Schr\"oder numbers while without it we would
 have the {\em large Schr\"oder numbers}.

For $k=3$, we obtain the following matrix identity on Motzkin
numbers:
\begin{center}
\begin{equation} \label{3p}
\begin{bmatrix} 1 \\1 & 1\\2 & 2 & 1\\4 & 5 & 3 & 1\\9 & 12 & 9 & 4 & 1\\& & \cdots &&&\ddots   \end{bmatrix}
\begin{bmatrix} 1 \\ 2 \\ 3 \\ 4 \\ 5 \\ \vdots
\end{bmatrix} = \begin{bmatrix} 1 \\ 3 \\ 3^2 \\ 3^3 \\ 3^4 \\ \vdots
\end{bmatrix},
\end{equation}
\end{center}
where the first column is the sequence of Motzkin numbers, and
matrix $A=(a_{ij})$ is generated by the following recurrence
relation:
\[ a_{i,j}= a_{i-1, j-1} + a_{i-1, j} + a_{i-1, j+1}. \]
For $k=5$, we find the following matrix identity
\begin{equation}\label{5p}
\begin{bmatrix} 1 \\3 & 1\\10 & 6 & 1\\36 & 29 & 9 & 1\\137 & 132 & 57 & 12 & 1\\& & \cdots &&&\ddots  \end{bmatrix}
\begin{bmatrix} 1 \\ 2 \\ 3 \\ 4 \\ 5 \\ \vdots
\end{bmatrix} = \begin{bmatrix} 1 \\ 5 \\ 5^2 \\ 5^3 \\ 5^4 \\ \vdots
\end{bmatrix},
\end{equation}
where the first column  sequence A002212 in \cite{sloane}, which
has two interpretations, the number of $3$-Motzkin paths or the
number of ways to assemble benzene rings into a tree \cite{hr}.
Recall that a $3$-Motzkin path is a lattice path from $(0,0)$ to
$(n-1,0)$ that does not go below the $x$-axis and consists of up
steps $U=(1,1)$, down steps $D=(1,-1)$, and three types of
horizontal steps $H=(1,0)$. The above matrix $A=(a_{i,j})$ is
generated by the first column and the following recurrence
relation
\[ a_{i,j}= a_{i-1,j-1} + 3a_{i-1,j} +
a_{i-1,j+1}.\]
 We may prove  the above identities (\ref{3p}) and (\ref{5p}) by using
method of Riordan arrays.  So the natural question is to find a
matrix identity for the sequence $(1, k, k^2, k^3, \ldots)$. We
need the combinatorial interpretation of the entries in the matrix
in terms of weighted partial Motzkin paths, as given by Cameron
and Nkwanta \cite{cn}. To be precise, a partial Motzkin path, also
called a Motzkin path from $(0,0)$ to $(n,k)$ in \cite{cn}, is
just a Motzkin path but without the requirement of ending on the
x-axis.  A weighted partial Motkzin a partial Motzkin path with
the weight assignment  that  the horizontal steps are endowed with
a weight $k$ and the down steps are endowed with a weight $t$,
where $k$ and $t$ are regarded as positive integers. In this
sense, our weighted Motzkin paths can be regarded as further
generalization of $k$-Motzkin paths in the sense of $2$-Motzkin
paths and $3$-Motkzin paths \cite{BdLPP, deutschs, sloane}.

We introduce the notion of weighted free Motzkin paths which is a
lattice path consisting of Motzkin steps without the restrictions
that it has to end with a point on the $x$-axis and  it does not
go below the $x$-axis. We then give a bijection between weighted
free Motzkin paths and weighted partial Motzkin paths with an
elevation line, which leads to a matrix identity involving the
number of weighted partial Motzkin paths and the sequence $(1, k,
k^2, \ldots)$.  The idea of the elevation operation is also used
by Cameron and Nkwanta in their combinatorial proof of the
identity (\ref{eq2.2}) in a more restricted form. By extending our
argument to weighted partial Motzkin paths with multiple elevation
lines, we obtain a combinatorial proof of an identity recently
derived by Cameron and Nkwanta, in answer to their question.

We also give a generalization of the matrix identity (\ref{eq2.2})
and give a combinatorial proof by using colored Dyck paths.

\section{Riordan Arrays}

In this section, we give a brief introduction to the notion of
Riordan arrays \cite{shapiroba, shapirotr, sprugnoli}. Let us use
(\ref{eq2.2}) and (\ref{eq1.1}) as examples. Start with two
generating functions $g(x)=1+g_1x+g_2x^2+\cdots$ and
$f(x)=f_1x+f_2x^2+\cdots$ with $f_1\neq 0$. Let
$H=(h_{i,j})_{n,j\geq 0}$ be the infinite lower triangular matrix
with nonzero entries on the main diagonal, where
$h_{i,j}=[x^i](g(x)(f(x)^j)$ for $i\geq j$, namely, $h_{i,j}$
equals the coefficient of $x^i$ in the expansion of the series
$f(x)(g(x)^j)$. If an infinite lower triangular matrix $H$ can be
constructed in this way from two generating functions $g(x)$ and
$f(x)$, then it is called a {\em Riordan array} and is denoted by
$H=(g(x),f(x))=(g,f)$.

Suppose we multiply the matrix $H=(g,f)$ by a column vector
$(a_0,a_1,\cdots)^T$  and get a column vector
$(b_0,b_1,\cdots)^T$.  Let  $A(x)$ and $B(x)$  be the
generating functions for the sequences $(a_0,a_1,\cdots)$  and
$(b_0,b_1,\cdots)$ respectively.  Then the method of Riordan arrays
asserts that
\[
B(x)=g(x)A(f(x)).\]

For the matrix identity (\ref{eq2.2}), let $g(x)$ be the
generating function for Catalan numbers $(1, 2, 5, 14, \ldots)$:
\[ g(x) = {1 -2x-
\sqrt{1-4x} \over 2x^2}.\] Let $f(x)=xg(x)$. From the recurrence
relation (\ref{r2.2}) one may derive that the generating function
for the sequence in the $j$-th $(j\geq 1)$ column in the matrix in
(\ref{eq2.2}) equals $g(xg)^{j-1}$. Let $H$ be the Riordan array
$(g, xg)$. Since the generating function of $(1,2,3,4\cdots)^T$
equals $A(x)=\frac{1}{(1-x)^2}$, it follows that
$B(x)=g(x)A(xg(x))=\frac{1}{1-4x}$ equals the generating function
for the right hand side of (\ref{eq2.2}). Thus we obtain the
identity (\ref{eq2.2}).

Let us consider the matrix identity (\ref{eq1.1}). Let $g(x)$ be
the generating function for the little Schr\"oder numbers as given
by
\begin{equation} \label{sg}
g(x)={1-3x-\sqrt{1-6x+x^2}\over 4x^2},
\end{equation}
and let $f(x) =xg(x)$. Note that the generating function for the
sequence  $(1, 3, 7, 15, \ldots)$ equals
$A(x)=\frac{1}{(1-x)(1-2x)}$. From the recurrence relation
(\ref{r1.1}) one may verify that the matrix in (\ref{eq1.1}) is
indeed the Riordan array $(g, xg)$. Therefore, the generating
function for the right hand side of (\ref{eq1.1}) equals $
g(x)A(xg(x))=\frac{1}{1-6x}$, which implies (\ref{eq1.1}).

Using the same method, we can verify the matrix identity
(\ref{3p}) and (\ref{5p}). Since we are going to establish a
general bijection for weighted Motzkin numbers, here we omit the
proofs.

\section{Dyck path interpretation of (\ref{eq2.2})}

In this section, we present a combinatorial interpretation of the
matrix identity (\ref{eq2.2}) by using Dyck paths.  A {\em Dyck
path} of length $2n$ is a path going from the origin $(0,0)$ to
$(2n,0)$ using steps $U=(1,1)$ and down steps $D=(1,-1)$ such that
no steps is below the $x$-axis \cite{De, stanley}. The number of
Dyck paths of length $2n$ equals the Catalan number $C_n$.

For a Dyck path $P$, the points on the $x$-axis except for the
initial point are called return points. In this sense, the ending
point is always a return point. Formally speaking, a {\em
composition} of a Dyck path $P$ is sequence of Dyck path $(P_1,
P_2, \ldots, P_j)$ such that $P=P_1 P_2 \cdots P_j$, where $P_1,
P_2, \ldots, P_j$ are Dyck paths.  For a composition $(P_1, P_2,
\ldots, P_j)$ of a Dyck path $P$, its length is meant to be the
length of $P$ and $j$ is called the number of segments. We may
choose certain return points to cut off the Dyck paths into a
composition. We use the convention that the ending point is always
a cut point. Clearly, a Dyck path with one segment is an ordinary
Dyck path.

\begin{lemma}\label{dotheom}
 For $j\geq 2$,   we have the following recurrence relation
\begin{equation}\label{eq.1}
 d_{i,j}=d_{i-1,j-1} +
2d_{i-1,j} + d_{i-1,j+1}.
\end{equation}
\end{lemma}

\pf Let   $(P_1, P_2, \ldots, P_{j})$ be composition of a Dyck
path $P$ of length $2i$. Consider the following cases for $P_1$.
Case 1: $P_1=UD$. Then we get a composition length $2(i-1)$ with
$j-1$ segments: $(P_2, \ldots, P_j)$. Case 2: $P_1=QUD$ and $Q$ is
not empty. Then we get a composition $(Q, P_2, \ldots, P_j)$ of
length $2(i-1)$ and $j$ segments. Case 3: $P_1=U Q D$, $Q$ is not
empty. We get a composition $(Q,  P_2, \ldots, P_j)$ of length
$2(i-1)$ with $j$ segments. Case 4: $P_1=Q_1UQ_2D$, where $Q_1$
and $Q_2$ are not empty. Then we get a composition $(Q_1, Q_2,
P_2, \ldots, P_j)$ of length $2(i-1)$ with $j+1$ segments. Adding
up the terms in the above cases, we obtain the desired recursion
(\ref{eq.1}). \qed

From Lemma \ref{dotheom} one sees  that
 the entry $a_{i,j}$ in the
triangular matrix of the identity (\ref{eq2.2}) can be explained
as the number compositions of Dyck paths of length $2i$ that
contain $j$ segments. We remark that this combinatorial
interpretation can also be derived from the generating function of
the entries in the $j$-th column in of the matrix in
(\ref{eq2.2}). The following formula for $a_{i,j}$ has been
derived by Cameron and Nkwanta \cite{cn}:
\[ a_{i,j}= {j \over i} \, {2i \choose i-j} .\]

Let us rewrite the  matrix identity (\ref{eq2.2}) as follows
\begin{equation}\label{eq.2}
\sum_{j= 1}^i ja_{i,j}=4^{i-1}. \end{equation} A combinatorial
formulation of the above identity is given by Callan
\cite{callan}.

 We are now ready to give a combinatorial proof of
the above identity. Clearly, $4^n$ is the number of sequences of
length $n$ on four letters, say, $\{1, 2, 3, 4\}$. The term
$ja_{i,j}$ suggests that we should specify a segment in a
composition as a distinguished segment. We may use a star $*$ to
mark the distinguished segment. We call a composition with a
distinguished segment a {\em rooted} composition of a Dyck path.
Then $ja_{i,j}$ equals the number of rooted compositions of Dyck
paths of length $2i$ that contain $j$ segments.

 \begin{theorem}
There is a bijection $\phi$  between the set of rooted
compositions of Dyck paths of length $2i$ and the set of sequences
of length $i-1$ on four letters.
\end{theorem}

\pf Let $(P_1, P_2, \ldots, P_j)$ be a rooted composition of a
Dyck path $P$ of length $2i$. We proceed to construct a sequence
of length $i-1$ on the elements $\{ 1, 2,3, 4\}$. We now
recursively define a map $\phi$ from  rooted compositions of a
Dyck path $P$ of length $2i$ to sequences of length $i-1$ on $\{
1, 2, 3, 4\}$. For $i=1$, $P$ is unique, and the sequence is set
to be the empty sequence. We now assume that $i>1$. Let $(P_1,
\ldots, P_t^*, \ldots, P_j)$ be a rooted composition of $P$ with
$P_t^*$ being the distinguished segment.

We have the following cases.
 \begin{enumerate}

 \item[1.]  $P_1=UD$ and $t=1$. Then we set
 $\phi(P)=1\, \phi(P_2^*, P_3, \ldots, P_j)$.

 \item[2.] $P_1=UD$ and $t\not=1$. Then we set
 $\phi(P)=2\, \phi(P_2, \ldots, P_t^*\ldots, P_j)$.

 \item[3.]  $P_1=QUD$ and $Q$ is  a nonempty Dyck path.
 Set $\phi(P)=3\,
 \phi(Q^*, P_2,  \ldots, P_j)$ if $t=1$ and set
 $\phi(P_1, \ldots, P_j) =3 \, \phi (Q, P_2, \ldots, P_t^*, \ldots, P_j)$ if $t>1$.

 \item[4.] $P_1=Q_1UQ_2D$, where $Q_1$ and $Q_2$ are nonempty Dyck paths.
 Then we set
\[\phi(P_1, \ldots, P_j)=1\, \phi
 (Q_1, Q_2, P_2, \ldots, P_t^*, \ldots, P_j)
 \text{ if }
 t>1 ,\]
\[\phi(P)=1\, \phi
 (Q_1, Q_2^*, P_2, \ldots, P_j) \text{ if } t=1.\]

\item[5.] If $P_1=UQD$ and $Q$ is a nonempty Dyck path. Then we
set
\[\phi(P)=4\, \phi(Q, P_2, \ldots, P_t^*, \ldots,
P_j).\]
 \end{enumerate}

 In order to show that $\phi$ is a bijection, we  construct the
 reverse map of $\phi$.
 Let $w=w_1w_2\cdots w_{i-1}$ be a sequence of length $i-1$ on $\{1, 2, 3, 4\}$.
 If $i=1$, then it
 corresponds to $UD$. We now assume that $i>1$. Suppose that $w_2w_3\cdots w_{i-1}$
 corresponds to a rooted composition $(R_1, R_2, \cdots, R_m)$
 of a Dyck path $P$ of length $2(i-1)$ with $R_k$ being the distinguished segment.
We proceed to find a rooted composition $(P_1, P_2, \ldots, P_j)$
with $P_t$ being the distinguished segment such that $\phi(P_1,
P_2, \ldots, P_j) = w_1 \phi(R_1, R_2, \ldots, R_m)$.

    If $w_1=2$, we have $P_1=UD$ and $(P_2, P_3, \ldots, P_j) =
                                       (R_1, R_2, \ldots, R_m)$.
                      It follows that $t=k+1$ and $j=m+1$. Thus we
                      can recover $(P_1, P_2, \ldots, P_j)$.
     For the case $w_1=3$, we have $P_1=R_1UD$ and  $t=k$.
      Also, we can recover $(P_2, \cdots, P_j)$ from $(R_2,  \ldots, R_m)$.
For the case $w_1=4$, we  have $t=k$, $P_1=UR_1D$, and $(P_2,
\ldots, P_j)=(R_2, \ldots, R_m)$.

It remains to deal with the situation $w_1=1$, which involves
Cases 1 and 4 of the bijection. If $k=1$, then we have $t=1$,
$P_1=UD$ and $(P_2,  \ldots, P_j)=(R_1,\ldots, R_m)$. If $k=2$,
then we have $t=1$, $P_1=R_1UR_2D$ and $(P_2,  \ldots, P_j)= (R_3,
 \ldots, R_m)$. If $k>2$, then we have $t=k-1$,
$P_1=R_1UR_2D$, and $(P_2,  \ldots, P_j)=(R_3,  \ldots, R_m)$.
 Thus, we have shown that $\phi$ is a
bijection. \qed

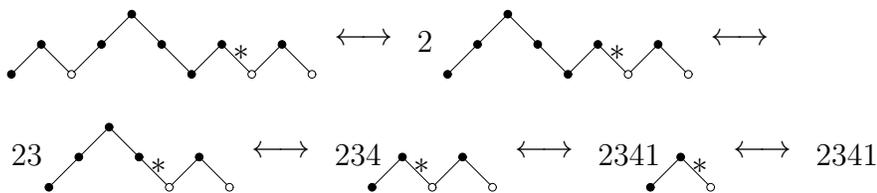
\begin{figure}[h,t]
\begin{center}
\setlength{\unitlength}{1mm}
\begin{picture}(120,25)

\put(0,15){\line(1,1){4}} \put(4,19){\line(1,-1){3.7}}
\put(8.3,15.3){\line(1,1){8}} \put(16,23){\line(1,-1){8}}
\put(24,15){\line(1,1){4}}
\put(28,19){\line(1,-1){3.7}}\put(29.5,16.9){$\ast$}\put(32.3,15.3){\line(1,1){4}}
\put(36,19){\line(1,-1){3.7}} \put(0,15){\circle*{1}}
\put(4,19){\circle*{1}} \put(8,15){\circle{1}}
\put(12,19){\circle*{1}} \put(16,23){\circle*{1}}
\put(20,19){\circle*{1}}
\put(24,15){\circle*{1}}\put(28,19){\circle*{1}}
\put(32,15){\circle{1}} \put(36,19){\circle*{1}}
\put(40,15){\circle{1}}

\put(43,19){$\longleftrightarrow$} \put(54,18){$2$}

\put(58,15){\line(1,1){8}} \put(66,23){\line(1,-1){8}}
\put(74,15){\line(1,1){4}}
\put(78,19){\line(1,-1){3.7}}\put(79.5,16.9){$\ast$}\put(82.3,15.3){\line(1,1){4}}
\put(86,19){\line(1,-1){3.7}} \put(58,15){\circle*{1}}
\put(62,19){\circle*{1}} \put(66,23){\circle*{1}}
\put(70,19){\circle*{1}}
\put(74,15){\circle*{1}}\put(78,19){\circle*{1}}
\put(82,15){\circle{1}} \put(86,19){\circle*{1}}
\put(90,15){\circle{1}}

\put(93,19){$\longleftrightarrow$}

\put(0,3){$23$} \put(5.3,0.3){\line(1,1){8}}
\put(13,8){\line(1,-1){7.7}}
\put(18.5,1.9){$\ast$}\put(21.3,0.3){\line(1,1){4}}
\put(25,4){\line(1,-1){3.7}} \put(5,0){\circle*{1}}
\put(9,4){\circle*{1}} \put(13,8){\circle*{1}}
\put(17,4){\circle*{1}} \put(21,0){\circle{1}}
\put(25,4){\circle*{1}} \put(29,0){\circle{1}}

\put(32,4){$\longleftrightarrow$}

\put(43,3){$234$} \put(48,0){\line(1,1){4}}
\put(52,4){\line(1,-1){3.7}}
\put(53.5,1.9){$\ast$}\put(56.3,0.3){\line(1,1){4}}
\put(60,4){\line(1,-1){3.7}} \put(48,0){\circle*{1}}
\put(52,4){\circle*{1}} \put(56,0){\circle{1}}
\put(60,4){\circle*{1}} \put(64,0){\circle{1}}

\put(67,4){$\longleftrightarrow$}

\put(78,3){$2341$} \put(85,0){\line(1,1){4}}
\put(89,4){\line(1,-1){3.7}} \put(90.5,1.9){$\ast$}
\put(85,0){\circle*{1}} \put(89,4){\circle*{1}}
\put(93,0){\circle{1}}

\put(96,4){$\longleftrightarrow$} \put(107,3){$2341$}

\end{picture}
\end{center}
\caption{The bijection $\phi$.} \label{dyck}

\end{figure}

An example of the bijection $\phi$ is given in Figure \ref{dyck},
where normal vertices are drawn with black points, return points
that cut off the Dyck paths into segments are drawn with white
points, and the distinguished segment is marked with a $*$ on its
last return step.

\section{Weighted Partial Motzkin Paths }

A {\em Motzkin path} of length $n$ is a path going from $(0,0)$ to
$(n,0)$ consisting of up steps $U=(1,1)$,  down steps $D=(1,-1)$
and horizontal steps $H=(1,0)$,  which never goes below the
$x$-axis. A {\em $(k,t)$-Motzkin path} is a
 Motzkin path such that each horizontal step is weighted  by $k$,
 each
down step is weighted by  $t$ and each up step is weighted by $1$.
The case $k=2, t=1$ gives the
 $2$-Motzkin paths which have been introduced  by
 Barcucci, del Lungo,  Pergola and  Pinzani \cite{BdLPP}
and have been studied by
 Deutsch and Shapiro  \cite{deutschs}.
 The {\em weight } of a  path is the product
 of the weights of all its steps. Denote by $|P|$ the weight of a path $P$.
  The {\em weight } of a set of  paths is the sum of the total
  weights of all the paths.
 For any step, we say that it is at level $k$ if the $y$-coordinate of
 its end point is at level $k$.

In this section, we aim to establish the following matrix identity
on weighted Motzkin numbers.

\begin{theorem} Let $M=(m_{i,j})_{i,j\geq 1}$ be the lower triangular matrix
such that the first column is the sequence of the number of
$(k-t-1, t)$-Motzkin paths of length $n$ and $m_{i,j}$ satisfies
the following recurrence relation for $j\geq 2$:
\begin{equation}\label{rc.5}
m_{i,j}=m_{i-1,j-1}+(k-t-1)m_{i-1,j}+tm_{i-1,j+1}.
\end{equation}
 Then we have
\begin{equation}\label{eq.5}
(m_{i,j})  \times
\begin{bmatrix} 1 \\ 1+t \\ 1+t+t^2  \\ \vdots\\
\end{bmatrix} = \begin{bmatrix} 1 \\ k \\ k^2  \\ \vdots
\end{bmatrix},
\end{equation}
\end{theorem}

It is well known that  the number of $2$-Motzkin paths of length
$n$ is given by the Catalan number $C_{n+1}$. It follows that the
matrix identity (\ref{eq2.2}) is a special case of  (\ref{eq.5})
for $k=4,t=1$. Denote by $f(x)=\sum_{n\geq 0}f_{n}x^n$ the
generating function for the number of $(3,2)$-Motzkin paths. Then
it is easy to find the functional equation for $f(x)$:
$$
f(x)=1+3xf(x)+2x^2f^2(x).
$$
It follows that
$$f(x)={{1-3x-\sqrt{1-6x+x^2}}\over {4x^2}}. $$
From the above generating function, one see that the number of
$(3,2)$-Motzkin paths of length $n$ equals the $n$-th little
Schr\"{o}der number.  Therefore, the matrix identity (\ref{eq1.1})
is a special case of (\ref{eq.5}) for $k=6,t=2$.

Let us rewrite matrix identity (\ref{eq.5}) in the following form
\begin{equation} \label{mi}
\sum_{j= 1}^i \, m_{i,j}(1+t+\cdots +t^{j-1}) =k^{i-1}.
\end{equation}

The following combinatorial interpretation of the entries in the
matrix in (\ref{eq.5}) is due to Cameron and Nkwanta \cite{cn}. A
partial $(k,t)$-Motkzin path is defined as an initial segment of a
$(k,t)$-Motkzin path. We say that a partial $(k, t)$-Motzkin path
ends at level $j$ if its
 last step is at level $j$.

\begin{lemma}[\cite{cn}] Let $m_{i,j}$ be the entries in the matrix in
(\ref{eq.5}). Then $m_{i,j}$ equals the the number of partial
$(k-t-1, t)$-Motzkin paths of length $i-1$ that end at level
$j-1$.
\end{lemma}

\noindent {\em Proof.} Regarding the first column of the matrix
$M$, one sees that a partial $(k-t-1, t)$-Motzkin path that ends
at level zero is just a $(k-t-1, t)$-Motzkin path. Let $a_{i,j}$
denote the number of partial $(k-t-1, t)$-Motzkin paths of length
$i-1$ ending at level $j-1$. Let $P$ be a partial $(k-t-1,
t)$-Motzkin path of length $i-1$ that ends at level $j-1$ $(j>1)$.
By considering the last step of $P$ and its weight, one sees that
$a_{i,j}$ satisfies the recurrence relation (\ref{rc.5}). \qed

Let $P$ be a partial $(k-t-1, t)$-Motzkin path  ending at level
$j-1$. We need the notion of an {\it elevated } partial Motzkin
path, which has been introduced by Cameron and Nkwanta \cite{cn}
in their combinatorial proof of the following identity which is a
reformulation of (\ref{eq.2}):
\[ 4^n=
\sum_{k=0}^n \, {(k+1)^2 \over n+1} \, {2n+2\choose n-k}.\]
 Let $p$  be an
integer  with $0\leq p\leq j$. The {\em elevation} of $P$ with
respect to the horizontal line $y=p$ is defined as follows. For
$p=0$, the elevation of $P$ with respect to $y=0$ is just $P$
itself. We now assume $ 0<p\leq j$. Note that there are always up
steps of $P$ at levels $j-1$, $j-2$, \ldots, $0$ bearing in mind
that an up step is said to be at level $k$ if its initial point is
at level $k$. Therefore, for each level from $0$ to $p-1$, one can
find a rightmost up step. Note that there are no other steps at
the same level to the right of the rightmost up step which is
called a {\em R-visible} up step with respect to the line $y=p$ in
the sense that it can seen far away from the right. The elevation
of $P$ with respect to the line $y=p$ is derived from $P$ by
changing the R-visible up steps up to level $p-1$ to down steps by
elevating their initial points. The line $y=p$ is called an {\em
elevation line}.

Figure \ref{mot} is an illustration of the elevation of a partial
Motzkin path with respect to the line $y=2$.

\begin{figure}
\begin{center}
\setlength{\unitlength}{0.9mm}
\begin{picture}(150,30)
\thicklines \put(-0.15,0){\line(1,1){4}}
\put(0.15,0){\line(1,1){4}}\thinlines \put(4,4){\line(1,0){4}}
 \put(8,4){\line(1,1){4}}
\put(12,8){\line(1,-1){4}} \put(16,4){\line(1,0){4}}
\put(20,4){\line(1,1){8}} \put(28,12){\line(1,-1){8}} \thicklines
\put(35.85,4){\line(1,1){4}}
\put(36.15,4){\line(1,1){4}}\thinlines
\put(40,8){\line(1,1){8}}\put(48,16){\line(1,-1){4}}
\put(52,12){\line(1,0){4}} \multiput(0,8)(3,0){20}{\line(1,0){1}}

\put(0,0){\circle*{0.8}} \put(4,4){\circle*{0.8}}
\put(8,4){\circle*{0.8}} \put(12,8){\circle*{0.8}}
\put(16,4){\circle*{0.8}} \put(20,4){\circle*{0.8}}
\put(24,8){\circle*{0.8}} \put(28,12){\circle*{0.8}}
\put(32,8){\circle*{0.8}} \put(36,4){\circle*{0.8}}
\put(40,8){\circle*{0.8}} \put(44,12){\circle*{0.8}}
\put(48,16){\circle*{0.8}} \put(52,12){\circle*{0.8}}
\put(56,12){\circle*{0.8}}

\put(63,5){$\longleftrightarrow$}

\thicklines \put(79.85,8){\line(1,-1){4}}
\put(80.15,8){\line(1,-1){4}}\thinlines \put(84,4){\line(1,0){4}}
 \put(88,4){\line(1,1){4}}
\put(92,8){\line(1,-1){4}} \put(96,4){\line(1,0){4}}
\put(100,4){\line(1,1){8}} \put(108,12){\line(1,-1){8}}
\thicklines \put(115.85,4){\line(1,-1){4}}
\put(116.15,4){\line(1,-1){4}}\thinlines
\put(120,0){\line(1,1){8}} \put(128,8){\line(1,-1){4}}
\put(132,4){\line(1,0){4}}

\put(80,8){\circle*{0.8}} \put(84,4){\circle*{0.8}}
\put(88,4){\circle*{0.8}} \put(92,8){\circle*{0.8}}
\put(96,4){\circle*{0.8}} \put(100,4){\circle*{0.8}}
\put(104,8){\circle*{0.8}} \put(108,12){\circle*{0.8}}
\put(112,8){\circle*{0.8}} \put(116,4){\circle*{0.8}}
\put(120,0){\circle*{0.8}} \put(124,4){\circle*{0.8}}
\put(128,8){\circle*{0.8}} \put(132,4){\circle*{0.8}}
\put(136,4){\circle*{0.8}}

\end{picture}
\end{center}
\caption{ Elevation of a partial Motzkin path}\label{mot}
\end{figure}
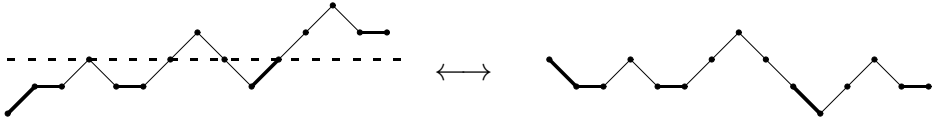

We now introduce the notion of {\em free Motzkin } paths which are
lattice paths starting from $(0,0)$ and using up steps $U=(1,1)$,
down steps $D=(1,-1)$ and  horizontal steps $H=(1,0)$. Note that
there is no restriction so that the paths may go below the
$x$-axis. A free $(k,t)$-Motzkin path is a free Motzkin path in
which the steps are weighted in the same way as for
$(k,t)$-Motzkin paths, namely, an up step has weight one, a
horizontal step has weight $k$ and a down step has weight $t$.

For a free Motzkin path $P$ we may analogously define the {\em
L-visible} down steps as the down steps that are visible from the
far left. It is clear that a complete Motzkin path (a partial
Motzkin path with ending point on the $x$-axis) has no R-visible
up steps. Similarly, a partial Motzkin path has no  L-visible down
steps.

We have the following summation formula for weighted free Motzkin
paths.

\begin{lemma}
The sum of weights of  free $(k-t-1, t)$-Motkzin paths of length
$i$ equals $k^i$.
\end{lemma}

The proof of the above lemma is obvious because of the relation
\[ (1+k-t-1+t)^i=k^i.\]
 We are now led to establish a bijection for
the identity (\ref{mi}).

\begin{theorem}\label{the.1}
There is a bijection between the set of partial
$(k-t-1,t)$-Motzkin paths of length $i$  with an elevation line
and the set of free $(k-t-1,t)$-Motzkin path of length $i$.
\end{theorem}

\pf The bijection is just the elevation operation. The reverse map
is also easy. For a free Motzkin path, one can identify the
L-visible down steps, if any, then change these L-visible down
steps to up steps by elevating their end points. \qed

For a partial $(k-t-1, t)$-Motzkin path $P$ with an elevation line
$y=p$, suppose that $Q$ is the elevation of $P$ with respect to
$y=p$. It is clear that the weight of $Q$ equals $t^p |P|$. If $P$
ends at level $j$, then the possible elevation lines are $y=0$,
$y=1, \ldots, y=j-1$. Summing over $j$ to get what we have
currently. Thus we arrive at a combinatorial interpretation of the
identity (\ref{mi}).

 As a consequence of Theorem \ref{the.1}, we obtain
the matrix identity (\ref{eq.5}). Note that the matrix identity
(\ref{eq2.2}) is a special case of (\ref{eq.5}) for  $k=4$ and
$t=1$, and the   matrix identity (\ref{eq1.1}) is also a special
case of  (\ref{eq.5})  for  $k=6$ and $t=2$.

\section{An Identity of Cameron and Nkwanta}

In their study of involutions in Riordan groups,  Cameron and
Nkwanta \cite{callan} obtained the following identity for $m\geq
0$, and asked for a purely combinatorial proof:
\begin{equation*}
\binom{n}{m}4^{n-m}=\sum_{k=0}^{n}\frac{k+1}{n+1}\binom{k+m+1}{k-m}\binom{2n+2}{n-k}.
\end{equation*}

It is clear that  identity (\ref{eq.2}) is a special case for
$m=0$. To be consistent with our notation, we may rewrite the
above identity in the following form:
\begin{equation}\label{cam}
\binom{i-1}{m}4^{i-1-m}=\sum_{j=1}^{i}\frac{j}{i}\binom{j+m}{2m+1}\binom{2i}{i-j}.
\end{equation}

We now give a bijective proof (\ref{cam}).

 We recall that the number of
 partial $2$-Motzkin paths of length $i-1$  ending at level
 $j-1$ is given by $a_{i,j}=\frac{j}{i}\binom{2i}{i-j}$.
We now consider the set of partial $2$-Motzkin paths with $m$
marked R-visible up steps and $m+1$ elevation lines such that
there is exactly one marked R-visible up step between two adjacent
elevation lines. We now have a combinatorial interpretation of the
summand in (\ref{cam}).

\begin{lemma}
The summand in (\ref{cam}) counts partial $2$-Motzkin paths of
length $i-1$ ending at level
 $j-1$ with $m$ marked R-visible up steps and $m+1$ elevation lines such that there is
 exactly one marked step between two adjacent elevation lines.
\end{lemma}

\pf Let $P$ be a partial $2$-Motzkin path of length $i-1$ ending
at level $j-1$. Suppose that there are $m$ marked R-visible up
steps with initial points at levels $j_1, j_2, \ldots, j_m$. Let
$t_1=j_1$  and  $t_{i}=j_{i}-j_{i-1}-1$ with $j_{m+1}=j-1$ for
$i\geq 2$. Then one see that the number of ways to choose the
$m+1$ elevation lines such that there is exactly $m+1$ marked
R-visible up step is equals to
\[  (t_1+1) (t_2+1) \cdots (t_m+1). \]
Note that the $t_i$'s range over $t_1+t_2+\cdots +t_{m+1}=j-m-1$.
Thus, the number of partial $2$-Motzkin paths of length $i-1$
ending at level $j-1$ with the required marked steps and elevation
lines equal
\begin{equation}\label{c.1}
\sum_{t_1+t_2+\ldots t_{m+1}=j-1-m}(t_1+1)(t_2+1)\cdots
(t_{m+1}+1).
\end{equation}
Let $g(x)=\sum_{n\geq 0}(n+1)x^n$. It is clear that $g(x)={1\over
(1-x)^2}$. Hence the summation (\ref{c.1}) equals the coefficient
of $x^{j-1-m}$ in the expansion of ${1\over (1-x)^{2m+2}}$, that
is, the binomial coefficient ${j+m\choose 2m+1}$. \qed

In fact, we may give a combinatorial interpretation of the
binomial coefficient ${j+m \choose 2m+1}$ in the above proof. Let
$P$ be a partial $2$-Motzkin path of length $i-1$ ending at level
$j-1$ with $m$ marked R-visible up steps and $m+1$ elevation lines
such that there is exactly one marked up step between two adjacent
elevation lines. Suppose that the $k$-th elevation line
  and the $k$-th marked up step of $P$ are at  level $x_k$ and  $y_k$,
   respectively.
Such a configuration can be represented as follows:
\[   t_1\, | \, t_2 \, * \,  t_3\,  |\,  t_4\, * \,  t_5\, | \, \cdots \,  |
\,  t_{2m}  \,  * \, t_{2m+1} \,  |\,  t_{2m+2} ,\] where $t_i$
denotes the numbers of unmarked R-visible up steps.     It is
clear that we have $t_1+t_2+\ldots+t_{2m+2}=j-1-m$, and the number
of solutions of this equation equals the numbers of ways to
distribute $j-1-m$ balls into $2m+2$ boxes when a box can have
more than one ball. So this number equals the binomial coefficient
${ j+m \choose 2m+1}$.

 We are now ready to give a combinatorial proof of
the identity of Cameron and Nkwanta. We recall that a $2$-Motzkin
path have two kind of horizontal steps, straight steps and wavy
steps. We now need to introduce  the third kind of horizontal
steps -- dotted steps. Therefore, the left hand side of
(\ref{cam}) is the number of free $3$-Motzkin paths  with exactly
$m$ dotted horizontal steps. We now give the following bijection
that leads to a combinatorial interpretation of (\ref{cam}).

\begin{theorem}\label{cameron}
There is a bijection  between  partial $2$-Motzkin paths of length
$i$ with $m$ marked R-visible up steps and $m+1$ elevation lines
such that there is exactly one marked step between two adjacent
elevation lines and free $3$-motzkin paths of length $i$ with
exactly $m$ dotted horizontal steps.
\end{theorem}

\pf Suppose that $P=P_1U^*P_2U^*\ldots P_{m}U^*P_{m+1}$ is a
partial $2$-Motzkin path with $m$ marked R-visible up steps and
$m+1$ elevation lines, we get a free $3$-motzkin path by changing
all the marked up steps to dotted horizontal steps and applying
the elevation operation for each $P_k$.

  Conversely,  given a free $3$-Motzkin path
$P=P_1  \dashrightarrow P_2\dashrightarrow \cdots  P_{m}
\dashrightarrow P_{m+1}$ with $m$ dotted horizontal steps, where
$\dashrightarrow$ denotes a dotted horizontal step,  then we can
get a partial $2$-motzkin path by changing each dotted horizontal
step to a marked up step and the L-visible down steps of each
$P_k$ to up steps by elevating their end points. \qed

\begin{figure}
\begin{center}
\setlength{\unitlength}{0.9mm}
\begin{picture}(138,30)
\thicklines \put(-0.15,0){\line(1,1){4}}
\put(0.15,0){\line(1,1){4}}\thinlines \qbezier(4,4)(4.5,5)(5,4)
\qbezier(5,4)(5.5,3)(6,4) \qbezier(6,4)(6.5,5)(7,4)
\qbezier(7,4)(7.5,3)(8,4) \put(8,4){\line(1,1){4}}
\put(12,8){\line(1,0){4}} \put(16,8){\line(1,-1){4}}
\put(20,4){\line(1,1){8}} \put(28,12){\line(1,-1){8}}
\put(36,4){\line(1,1){4}} \thicklines \put(39.85,8){\line(1,1){8}}
\put(40.15,8){\line(1,1){8}}\thinlines
\put(42.5,9){\scriptsize$*$} \put(40,8){\line(1,1){12}}
\put(52,20){\line(1,0){4}} \multiput(0,4)(3,0){20}{\line(1,0){1}}
\multiput(0,16)(3,0){20}{\line(1,0){1}} \put(0,0){\circle*{0.8}}
\put(4,4){\circle*{0.8}} \put(8,4){\circle*{0.8}}
\put(12,8){\circle*{0.8}} \put(16,8){\circle*{0.8}}
\put(20,4){\circle*{0.8}} \put(24,8){\circle*{0.8}}
\put(28,12){\circle*{0.8}} \put(32,8){\circle*{0.8}}
\put(36,4){\circle*{0.8}} \put(40,8){\circle*{0.8}}
\put(44,12){\circle*{0.8}} \put(48,16){\circle*{0.8}}
\put(52,20){\circle*{0.8}} \put(56,20){\circle*{0.8}}
\put(63,5){$\longleftrightarrow$} \thicklines
\put(79.85,8){\line(1,-1){4}}
\put(80.15,8){\line(1,-1){4}}\thinlines
\qbezier(84,4)(84.5,5)(85,4) \qbezier(85,4)(85.5,3)(86,4)
\qbezier(86,4)(86.5,5)(87,4) \qbezier(87,4)(87.5,3)(88,4)
\put(88,4){\line(1,1){4}} \put(92,8){\line(1,0){4}}
\put(96,8){\line(1,-1){4}} \put(100,4){\line(1,1){8}}
\put(108,12){\line(1,-1){8}} \put(116,4){\line(1,1){4}}
\thicklines \multiput(120,8)(1,0){5}{\line(1,0){0.4}}
\put(123.85,8){\line(1,-1){4}}
\put(124.15,8){\line(1,-1){4}}\thinlines
\put(128,4){\line(1,1){4}} \put(132,8){\line(1,0){4}}
\put(80,8){\circle*{0.8}} \put(84,4){\circle*{0.8}}
\put(88,4){\circle*{0.8}} \put(92,8){\circle*{0.8}}
\put(96,8){\circle*{0.8}} \put(100,4){\circle*{0.8}}
\put(104,8){\circle*{0.8}} \put(108,12){\circle*{0.8}}
\put(112,8){\circle*{0.8}} \put(116,4){\circle*{0.8}}
\put(120,8){\circle*{0.8}} \put(124,8){\circle*{0.8}}
\put(128,4){\circle*{0.8}} \put(132,8){\circle*{0.8}}
\put(136,8){\circle*{0.8}}
\end{picture}
\end{center}
\caption{ Elevation of a partial $2$-Motzkin path}\label{mot}
\end{figure}
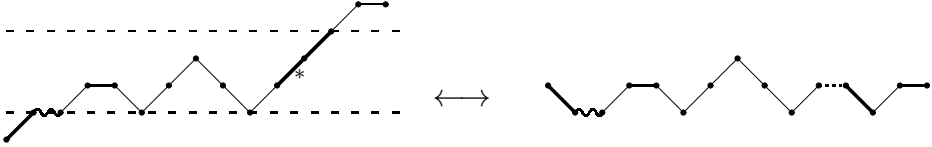

We conclude this section by giving a more general identity. Let
$a_{i,j,k}$ be the number of partial $k$-Motzkin paths of length
$i-1$  ending at level $j-1$. Then we have
\begin{equation}\label{cameron1}
\binom{i-1}{m}k^{i-1-m}=\sum_{j=1}^{i}a_{i,j,k-2}\binom{j+m}{2m+1}.
\end{equation}

\section{A Dyck path generalization of (\ref{eq1.1})}

In this section, we give a Dyck path generalization of the matrix
identity (\ref{eq1.1}) on the little Schr\"oder numbers. A
$k$-Dyck path is a Dyck path in which an up step is colored by one
of the $k$ colors $\{ 1, 2, \ldots, k\}$ if it not immediately
followed by a down step . In this section, we aim to give the
following generalization of (\ref{eq1.1}).

\begin{theorem}
Let  $M=(m_{i,j})$ be a lower triangular matrix with the first
column being the weight  of $(t^2-t)$-Dyck paths of length $2i$.
The other columns of $M$ are given by the following relation:
\begin{equation}\label{eq.6}
m_{i,j}=m_{i-1,j-1}+(t^2-t+1)m_{i-1,j}+(t^2-t)m_{i-1,j+1}.
\end{equation}
Then we have the following matrix identity
\begin{equation}\label{eq.t}
(m_{i,j}) \times
\begin{bmatrix} 1 \\ t^2-(t-1)^2 \\     \ t^3-(t-1)^3 \\ \vdots
\end{bmatrix} = \begin{bmatrix} 1 \\ t^2+t \\   (t^2+t)^{2} \\
\vdots
\end{bmatrix}.
\end{equation}
\end{theorem}

The matrix identity (\ref{eq1.1}) is a consequence of (\ref{eq.t})
by setting $t=2$. By using generating functions, one can verify
that the number of $2$-Dyck paths of length $2n$ equals the number
of little Schr\"oder paths of length $n$.

We now proceed to give a combinatorial proof of (\ref{eqt.t}). To
this end, we need to give a combinatorial interpretation of the
entries in the matrix $M$ in (\ref{eq.t}). We may define a
composition of $k$-Dyck path $P$ as a sequence of $k$-Dyck paths
$(P_1, P_2, \ldots, P_j)$ such that $P=P_1P_2\cdots P_j$, where
$j$ is called the number of segments.

\begin{lemma} Let $a_{i,j}$ be the sum of weights of
compositions of $(t^2-t)$-Dyck paths of length $2i$ with $j$
segments. Then $a_{i,j}$ satisfies the recurrence relation
(\ref{eq.6}).
\end{lemma}

The proof of the above lemma is similar to that of Lemma
\ref{dotheom}. Let us rewrite (\ref{eq.t}) as follows
\begin{equation}\label{eqt.t}
\sum_{j\geq 1}m_{i,j}(t^j-(t-1)^j)=(t^2+t)^{i-1}.
\end{equation}

In order to deal with $m_{i,j}(t^j-(t-1)^j)$ combinatorially, we
introduce a coloring scheme on a composition of a $(t^2-t)$-Dyck
path with $j$ segments. Suppose that we have $t$ colors $c_1, c_2,
\ldots, c_t$. If we use these $t$ colors to color the $j$ segments
such that the first color $c_1$ must be used, then there are
$t^j-(t-1)^j$ ways to accomplish such colorings.  We simply call
such colorings {\em $t$-feasible colorings}.

We can now present a bijection leading to a combinatorial proof of
(\ref{eqt.t}).

\begin{theorem}
There is a bijection between the set of compositions of
$(t^2-t)$-Dyck paths of length $2i$ with a $t$-feasible coloring
on the segments and the set of sequences of length $i-1$ on
$t^2+t$ letters.
\end{theorem}

The desired bijection $\sigma$ is constructed as follows. Let
$(P_1, P_2, \ldots, P_j)$ be a composition of a $(t^2-t)$-Dyck
path $P$ of length $2i$ with a $t$-feasible coloring on the
segments. We will use the following alphabet that contains $t^2+t$
letters: \begin{equation} \label{alphabet}
 \{ \alpha_r,  \;|\;
1\leq r\leq t\} \cup \{ \beta_s, \;|\; 1\leq s\leq t-1\} \cup \{
\gamma_{k} \; | \; 1 \leq k \leq t^2-t\} \cup \{ \delta\}.
\end{equation}

For $i=1$, both the composition and the $t$-feasible coloring are
unique. We set the corresponding sequence to be empty.
 For $i\geq 2$,  we consider the following cases:
 \begin{enumerate}
 \item[1.]If $P_1=UD$, $P_1$ is colored by $c_r$ ($1\leq r\leq t$)
 and $(P_2, \ldots, P_j)$
 still has a $t$-feasible coloring. Then we set
 $\sigma(P_1, \ldots P_j)=\alpha_r \sigma(P_2, \ldots, P_j)$.

 \item[2.]If $P_1=UD$, $P_1$ is colored by $c_1$ and $(P_2, \ldots, P_j)$
 does not inherit a
 $t$-feasible coloring. Assume that
 $P_2$ is colored by $c_{s+1}$ ($1\leq s\leq t-1$).
 Then we change the color of $P_2$ to $c_1$ and set
 $\sigma(P_1, \ldots, P_j)=\beta_s \sigma(P_2, \ldots, P_j)$.

 \item[3.]If $P_1=UDQ$, where $Q$ is not empty. Then we set
 $\sigma(P_1, \ldots, P_j)=\delta\sigma(Q, P_2, \ldots, P_j)$.

 \item[4.]If $P_1=UQD$,
 where $Q$ is not empty and the first  up step of $P$
 has color $k$ ($1 \leq k \leq
t^2-t$) because the first step of $Q$ is an up step. Then we set
$\sigma(P_1, \ldots, P_j)=\gamma_{k}
 \sigma(Q, P_2, \ldots, P_j)$.

\item[5.] If $P_1=U Q_1 D Q_2$, neither $Q_1$ nor $Q_2$ is empty
and the first up step of $P$ has color $c_k$. Since $k$ ranges
from $1$ to $t(t-1)$, we may encode a color $c_k$ by a pair of
colors $(c_p, c_q)$ where $p$ ranges from $1$ to $t$ and $q$
ranges from $1$ to $t-1$. Moreover, we may use $(c_r, \beta_s)$ to
denote a color $c_k$.  Then we assign color $c_r$
 to  $Q_1$, pass the color of $P_1$ to $Q_2$, and set
 $\sigma(P)=\beta_s\sigma(Q_1, Q_2, P_2, \ldots, P_j)$.
 \end{enumerate}
 For each case, the resulting path is always a sequence of length $i-1$.

In order to show that $\sigma$ is a bijection, we proceed to
construct the inverse map of $\sigma$. Let $S$ be a sequence of
length $i-1$ on the alphabet (\ref{alphabet}). If $i=1$, then we
get the unique Dyck path $UD$  and the unique composition with a
$t$-feasible coloring. Note that the up step in the Dyck path $UD$
is not colored. We now assume that $i>1$. It is easy to check that
 Cases 1, 3, and 4 are reversible. It remains to show that
Cases 2 and 5 are reversible. In fact, we only need to ensure that
Case 2 and Case 5 can be distinguished from each other. For Case
2, either $j=2$ or $(P_3, \ldots, P_j)$ does not have a
$t$-feasible coloring. On the other hand, for Case 5, $(Q_2,
P_2,\ldots, P_j)$ is always nonempty and it has a $t$-feasible
coloring. This completes the proof.  \qed

 We also have a combinatorial
interpretation of the matrix identity (\ref{eq1.1}) based on
little Schr\"{o}der paths. The idea is similar to the proof given
above, so the proof is omitted.

  \vskip 5mm

\noindent{\bf Acknowledgments.} This work was supported by the 973
Project on Mathematical Mechanization, the National Science
Foundation, the Ministry of Education, and the Ministry of Science
and Technology of China. The third author  is partially supported
by NSF grant HRD 0401697.



\begin{thebibliography}{100}

\bibitem{BdLPP}
E. Barcucci, A. del Lungo, E. Pergola and R. Pinzani, A
construction for enumerating $k$-coloured Motzkin paths, Lecture
Notes in Computer Science, Vol. 959, Springer, Berlin, 1995, pp.
254-263.

\bibitem{aigner}
M. Aigner, Catalan-like numbers and determinants, {\itshape J.
Combin. Theory,  Ser. A}, 87 (1999) 33-51.

\bibitem{callan} D. Callan,  A combinatorial interpretation of a
Catalan numbers identity, {\itshape Math. Mag.}, 72 (1999)
295-298.

\bibitem{cn}
N. Cameron and A. Nkwanta, On some (pseado) involutions in the
Riordan group, {\itshape J. Integer Sequences}, 8 (2005), Article
05.3.7.

\bibitem{chapman}
R. Chapman, Moments of Dyck paths, {\itshape Discrete Math.}, 204
(1999) 113-117.

\bibitem{De}
 E. Deutsch, Dyck path enumeration, {\itshape Discrete
Math.}, 204 (1999) 167-202.

\bibitem{deutschs}
E. Deutsch and L. Shapiro,  A bijection between ordered trees and
$2$-Motzkin paths and its many consequences, {\itshape Discrete
Math.}, 256 (2002) 655-670.

\bibitem{hr} F. Harary and R. C. Read,
            The enumeration of tree-like polyhexes, {\itshape
             Proc. Edinburgh Math. Soc.}, (2) 17 (1970) 1-13.

\bibitem{shapiroba}
L. Shapiro, Bijections and the Riordan group, {\itshape
Theoretical Computer Science}, 307 (2003) 403-413.

\bibitem{shapirotr}
L. Shapiro, S. Getu, Wen-Jin Woan, and L.C. Woodson, The Riordan
group, {\itshape Discrete Appl. Math.}, 34 (1991) 229-239.

\bibitem{shapirotc}
 L. Shapiro, Wen-Jin Woan, and  S. Getu, Runs, slides and moments, {\itshape
 SIAM J. Alg. Disc. Math.}, 4 (1983), 459-466.

\bibitem{shapiro}
 L. Shapiro,  A Catalan triangle, {\itshape Discrete Math.}, 14 (1976), 83-90.



\bibitem{sloane}
N.J.A. Sloane, S. Plouffe, The Encyclopedia of Integer Sequence,
Academic Press, San Diego, 1995, online at {\tt
www.research.att.com/\char126 njas/sequences/}.

\bibitem{sprugnoli}
R. Sprugnoli, Riordan arrays and combinatorial sums, {\itshape
Discrete Math.}, 132 (1994) 267-290.

\bibitem{stanley}
R.P. Stanley, {\itshape Enumerative Combinatorics, Vol. 2},
Cambridge University Press, Cambridge, 1999.

\bibitem{sulanke}
R.A. Sulanke, Moments of generalized Motzkin pahts, {\itshape J.
Integer Sequences}, 3 (2000) 00.1.1.



\bibitem{woan}
Wen-Jin Woan, L. Shapiro, and D.G. Rogers, The Catalan numbers,
the Lebesgue integral, and $4^{n-2}$, {\itshape The Amer. Math.
Monthly}, 104 (1997) 927-931.


\end{thebibliography}
\end{document}